\documentclass[12pt]{article}
\usepackage{epsfig,amsmath,amsfonts,amssymb,graphicx}

\newcommand{\beq}{\begin{quote}}
\newcommand{\enq}{\end{quote}}
\newcommand{\be}{\begin{equation}}
\newcommand{\en}{\end{equation}}
\newcommand{\del}{\delta}

\newcommand{\lag}{{\cal L}}
\newcommand{\x}{{\mathbf x}}
\newcommand{\F}{{\mathbf F}}
\newcommand{\J}{{\mathbf J}}
\newcommand{\p}{\partial}
\newcommand{\Z}{\Bbb Z}
\newcommand{\R}{\Bbb R}
\DeclareMathAlphabet{\varmathbb}{U}{bbold}{m}{n}
\newcommand{\one}{\varmathbb 1}

\newcommand{\deltau}{\delta\tau}
\newcommand{\dels}{\delta}

\newcommand{\mat}{\begin{pmatrix}}
\newcommand{\rix}{\end{pmatrix}}

\newcommand{\old}[1]{}

\begin{document}

\title{New Periodic Orbits for the $n$-Body Problem}
\author{Cristopher Moore\\
Computer Science Department \\
and Department of Physics and Astronomy\\
University of New Mexico, Albuquerque \\ 
and the Santa Fe Institute \\
\and Michael Nauenberg\\
Physics Department \\
University of California, Santa Cruz}
\maketitle

\begin{abstract}
Since the discovery of the figure-8 orbit for the three-body problem 
[Moore 1993] a large number of periodic orbits of 
the $n$-body problem with equal masses and beautiful symmetries have been discovered.  
However, most of those that have appeared in the literature 
are either planar or are obtained from  perturbations of planar orbits.  
Here we exhibit a number of new three-dimensional periodic 
$n$-body orbits with equal masses and cubic symmetry, 
including some whose moment of inertia tensor is a scalar.  
We found these orbits numerically by minimizing the action as 
a function of the trajectories' Fourier coefficients.  
We also give numerical evidence that a planar 3-body orbit first found 
in [H\'enon, 1976], rediscovered by [Moore 1993], and found to
exist for different masses by [Nauenberg 2001], is dynamically stable.

It is a pleasure to dedicate this paper to Philip Holmes.
\end{abstract}

\section{Introduction}

Inspired by the work of Holmes and his collaborators on the classification of periodic orbits in three-dimensional flows according to their knot and braid types~\cite{holmes1,holmes2,holmes3,holmesbook}, in 1993 Moore~\cite{moore} considered periodic orbits of the $n$-body problem in the plane, and asked which braids can appear in their space-time trajectories.  We can search for such orbits numerically by starting with a fictional orbit with a desired symmetry and topology, and then minimizing the action $S=\int \lag \,dt$ using gradient descent, where $\lag = K-V$ is the Lagrangian and $K$ and $V$ are the kinetic and potential energies respectively.  If $\tau$ parametrizes the progress of this descent, then with each step $\deltau$ we increment the position $\x_i(t)$ of the $i$th mass at each time $t$ by $\delta\x_i(t)$, where
\be
\delta\x_i(t) = - \deltau \frac{\delta S}{\delta \x_i(t)} = \deltau \bigl( m_i \ddot \x_i(t) - \F_i(t) \bigr)
\label{eq:descent}
\en
where $\F_i(t) = -\p V/\p \x_i(t)$ is the force acting on the $i$th mass 
at time $t$, and $\ddot \x_i(t) = d^2 \x_i/dt^2$ is its acceleration.  
Each time we update $\x_i(t)$ in this way, the action $S$ decreases.  This process can lead to one of three outcomes as $\tau$ increases towards infinity: escape, in which $\x_i(t)$ tends to infinity for some $i$ and some $t$; collision, in which $\x_i(t)=\x_j(t)$ for some $i \ne j$ at some $t$; or success, in which the $\x_i(t)$ 
converge to an orbit which is a (local) minimum of the action, and therefore a genuine trajectory 
of the equations of motion where $\F_i=m_i \ddot \x_i$ for all $i$.

For homogeneous potentials $V(r) \propto r^{\alpha}$, it is easy to show~\cite{moore} that escape cannot occur if $\alpha < 2$ unless the braid can be separated into two pieces,  
and collision cannot occur (starting with a fictional orbit with finite action) if $\alpha \le -2$ since the action of a colliding trajectory in that case is infinite.  Thus {\em all} braid types occur for $\alpha \le -2$, including the so-called ``strong force'' $\alpha = -2$, but only some subset of the possible braids will exist for, say, Newtonian gravity, where $\alpha = -1$.  

\begin{figure}
\centerline{
\includegraphics[width=5in]{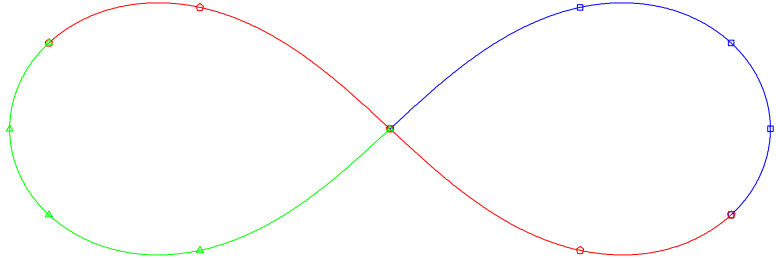}
}
\caption{The figure-8 orbit for three equal masses. The points show the masses at five equal time intervals.}
\label{fig:8}
\end{figure}

Performing this minimization on discretized orbits, Moore~\cite{moore} found a number of new periodic 3-body orbits in the plane, with a variety of braid types and symmetries: in particular, an orbit in which three equal masses chase each other around a figure-8 as shown in Figure~\ref{fig:8}, describing a classic three-strand braid in space-time.  Surprisingly, simulations indicated that this orbit is stable to small perturbations.  Six years later, this orbit was rediscovered by Chenciner and Montgomery~\cite{cm}, who provided a rigorous proof of its existence.  More precise numerical work was done by Sim\'o~\cite{simo1} and by Nauenberg~\cite{michael}; Sim\'o identified many new orbits, called ``choreographies,'' where masses follow each other along a fixed trajectory, and provided additional numerical evidence that the figure-8 orbit is stable to perturbations in the plane~\cite{simo2}. 

However, naively performing gradient descent in the variables $x_i(t)$ can introduce numerical problems.  In particular, discretizing the second time derivative as $\ddot \x(t) = \x(t+1) - 2 \x(t) + \x(t-1)$ produces a ``zig-zag'' instability in which the orbit jumps back and forth with alternate steps in $t$ and $\tau$.  Moreover, this method can only find local minima of the action, and thus is useless if the solution is another type of extremum.

An alternative approach was considered by Nauenberg~\cite{michael} which capitalizes on the periodicity of the coordinates $x_i(t)$, by expanding these variables in a Fourier series, and then finding extrema of the action $S$ as a function of the Fourier coefficients of the trajectory.  These extrema can be found by a gradient procedure which is analogous to the approach in~\eqref{eq:descent}, but avoids the instabilities associated with that method, and is not restricted to finding minima of $S$.  For forces associated with the homogeneous potential $V \propto r^\alpha$, the orbits scale as $x(t)=T^{2/(2-\alpha)}f(t/T)$ where $T$ is the period, and for our discussion we choose $T=2\pi$.  We substitute in $S$ the Fourier expansion 
\be 
\label{eq:fou}
x(t)=\sum_k a_k \sin kt
\en  
for each mass, and similarly for $y(t)$ and $z(t)$ (we suppress the subscript $i$ and ignore cosine terms for simplicity).  Then we have, to first order in a small change $\delta a_k$ in the Fourier coefficients $a_k$,  
\be
\delta S= \sum \frac{\partial S}{\partial a_k}\delta a_k,
\en
As in~\eqref{eq:descent} we perform gradient descent, but now with the rate of descent in each coefficient $a_k$ controlled by a parameter $\deltau_k$ which can depend on $k$:
\be
\label{eq:dela}
\delta a_k = - \deltau_k \frac{\p S}{\p a_k} \enspace . 
\en
In particular, substituting $K=(1/2)m(dx/dt)^2$ for the kinetic energy in $S$ yields
\be
\label{eq:delta1}
\delta a_k = -\deltau_k \left(m k^2 a_k +\int_0^{2\pi} dt \,\frac{\p V}{\p x} \sin kt \right) \enspace .
\en
Starting with given values of the Fourier coefficients $a_k$ and iterating with finite steps $\deltau_k$ leads to new values of $a_k$ which decrease the action $S$, leading as before to a local minimum of $S$.  Moreover, by setting some Fourier coefficients to zero or imposing equalities between the coefficients of different masses, we can restrict our search to orbits with certain symmetries, and by choosing the starting trajectory we can search for orbits with certain topologies.

We recover the naive gradient descent equation of~\cite{moore} by setting $\deltau_k = \deltau$ independent of $k$.  Note from~\eqref{eq:fou} that $\delta\x(t) = \sum_k \delta a_k \sin kt$; then substituting~\eqref{eq:delta1} yields~\eqref{eq:descent}.  However, allowing $\deltau_k$ to depend on $k$ allows us to avoid the above-mentioned zig-zag instability.  This instability, which is associated with discretizing the second time derivative of $\x$, is reflected here in the quadratic term in $k$ on the right hand side of~\eqref{eq:delta1}; if we set $\deltau_k=\dels/mk^2$,  where $\dels$ is a small positive constant, this removes this high-frequency instability, and gives a procedure which converges to a smooth trajectory whose Fourier coefficients $a_k$ tend to zero as $k \to \infty$.  In addition to removing the high-frequency instability, allowing $\deltau_k$ to depend on $k$ can, in some cases, allow us to find solutions which are extrema but not minima of $S$: if $\deltau_k$ is negative for some values of $k$, then we can sometimes find orbits which are saddle points of the action.  

Using this approach, Nauenberg~\cite{michael,michael2} found extensions in three dimensions of the figure-8 which have nonzero angular momentum about the two axes of symmetry in the original plane of the eight, and are periodic in a rotating frame; in particular, for the rotating-8 around the $y$ axis, these were found with negative values of $\deltau_k$ for $k=1$ and $k=3$.  Likewise, he found extended orbits with angular momentum normal to this plane, which remain planar.  For the case where the axis of rotation is along the horizontal symmetry axis of the figure-8, a proof for the existence of these orbits was given by C. Marchal~\cite{march1} who showed that these orbits interpolate between the figure-8 (which has zero total angular momentum) and the Lagrange orbit (where three masses co-rotate around a common center); see also~\cite{chenciner}.  

In this paper we present several new three-dimensional periodic orbits.  These orbits have symmetry that is fully three-dimensional, unlike, for instance, the ``hip-hop'' orbit of~\cite{hiphop}, which is obtained by perturbing a planar Lagrange orbit and has only dihedral symmetry.  In particular, we find a family of orbits with $4m$ masses where $m$ is an odd integer, with $m$ masses on each of four loops which are related by cubic symmetry. These orbits have zero total angular momentum.  Moreover, if $m$ is a multiple of $3$ they have cubic symmetry at all times, and their mass distribution is spherically symmetric to second order.  Although we have not attempted to do this, we note that it should be possible to prove the existence of these orbits rigorously using techniques like those in~\cite{march2,chen,ferrario}.

We then discuss another orbit we call the ``criss-cross,'' in which two masses orbit each other while the third rotates around them both in the opposite direction.  In fact, this orbit was first found by H\'enon as a member of a finite angular momemtum family of orbits periodic in rotating frames, extended from the zero angular momentum Schubart orbit~\cite{henon}.  It was rediscovered by Moore~\cite{moore} by searching for orbits with a particular braid type, and later extended to different masses by Nauenberg~\cite{michael}.  Recently, its existence was proved rigorously by Chen~\cite{chen} as part of a family of retrograde orbits.    Here we give numerical evidence that this orbit is dynamically stable in three dimensions, allowing the tantalizing possibility that it can be realized in astronomical systems.


\section{Orbits with cubic symmetry}

\begin{figure}
\centerline{
\includegraphics[width=5in]{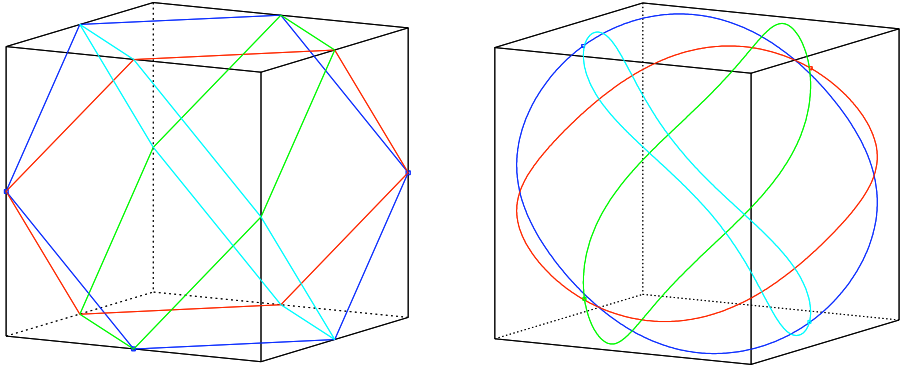}
}
\caption{Left, a cuboctahedron, and its four hexagonal ``great circles''; right, the 12-mass orbit with cubic symmetry.}
\label{fig:12}
\end{figure}

One source of inspiration for symmetric orbits in three dimensions consists of polyhedra.  For example, truncating the corners of a cube until a pair of triangular faces meet at the midpoint of each edge yields an Archimedean solid, the {\em cuboctahedron}, shown on the left in Figure~\ref{fig:12}.  Its 24 edges can be thought of as composing 4 hexagonal ``great circles'', and these can be oriented in such a way that half of the triangular faces are oriented counterclockwise and the others are oriented clockwise.  It takes only a little imagination to suppose that masses can start at the the corners of the cuboctahedron, and then travel along these great circle loops in these directions: first orbiting around the counterclockwise triangles, and then around the clockwise ones.  

The simplest such orbit has 4 masses, one on each of the four loops.  At all times, the masses are related to each other by $\pi$-rotations about the $x$, $y$ and $z$ axes: that is, by the transformations 
\be 
\label{eq:klein} (x,y,z) \to (x,y,z), (x,-y,-z), (-x,y,-z), (-x,-y,z) \enspace . 
\en
Each mass rotates counterclockwise around a diagonal axis, $(+1,+1,+1)$, $(+1,-1,-1)$, $(-1,+1,-1)$, or $(-1,-1,+1)$.  The total angular momentum is zero, a fact which we can immediately see from symmetry, since no nonzero vector $\J$ is fixed under the four transformations of~\eqref{eq:klein}. 

To find this orbit numerically, we construct circular orbits which satisfy these symmetry conditions, and use these as the starting point for the process of action minimization described in the Introduction.  For example, the coordinates of a circular orbit with unit radius normal to the diagonal $(1,1,1)$ are given by 
\be
x(t)=\sin t,\; y(t)=\sin (t+2 \pi/3),\; \mbox{ and } \; z(t)=\sin (t+4 \pi/3) \enspace . 
\en
Since our minimization process preserves the symmetries of the orbit, for all times $t$ the actual trajectory of this mass has the form
\be
x(t)=f(t),\; y(t)=f(t+2 \pi/3),\; \mbox{ and } \; z(t)=f(t+4 \pi/3) \enspace . 
\en
and the trajectories of the other masses are related to this one by the transformations~\eqref{eq:klein}.  

We now expand $f(t)$ in a Fourier series
\[ f(t) = \sum_k a_k \sin kt + b_k \cos kt \enspace , \]
and our goal is to find the Fourier coefficients $a_k$ and $b_k$.  However, if we set $f(0)=0$ so that the initial position of this mass is in the $y-z$ plane, by time reversal symmetry we have $f(t) = -f(-t)$ and the $b_k$ are identically zero.  Moreover, these orbits also have the symmetry that after half a period each mass is diametrically opposite to its previous position, so $f(t+\pi) = -f(t)$ and the Fourier coefficients $a_k$ are zero except for odd $k$.  

Minimizing the action, starting with these circular orbits and using this restricted set of Fourier coefficients, gives the periodic orbit whose Fourier coefficients are given in the first column of Table~\ref{tab:fourier}.  To verify adequate convergence of our Fourier sums, we compared the accelation of the resulting trajectory with the force for a series of values of $t$.  By considering frequencies up to $k \le 27$ we found agreement between the acceleration and the force up to $5$ decimal figures.  Higher accuracy can  be achieved by including higher frequencies in the Fourier sums.
We also integrated the equations of motion starting with the initial values we obtained, namely
\begin{equation}
\label{eq:cubic-initial}
 \x = (0, -0.69548, 0.69548), \; \dot \x = (0.87546, -0.31950, -0.31950) 
\end{equation}
using the Runge-Kutta algorithm.  We found good agreement during a single period; over longer times the integrated orbit rapidly diverges, indicating that this periodic orbit is unstable. 

Next, we consider orbits with $m$ masses on each of the four loops.  Now, in addition to the symmetries described above, these masses are related by the Lagrange symmetry, in which the $j$  masses are related by time translation, $t \to t + 2j \pi/m$ for $j=0,\ldots,m-1$.  

In fact, we will show that $m$ has to be an odd integer in order to avoid collisions between masses traveling on different loops.  By symmetry, each pair of loops intersects at two diametrically opposite points.  If $m$ is even, then the Lagrange symmetry implies that each mass has a ``partner'' at the diametrically opposite point on the same loop.  Hence, at the moment when one of the masses crosses an intersection with another loop, its partner crosses the opposite intersection.  However, the masses on different loops are related by the rotations~\eqref{eq:klein}, and so at the same moment this other loop also has a mass at each of these two intersections, and there will be two simultaneous collisions (and, by symmetry, collisions for other pairs of loops as well).

By the same reasoning, if there is an odd number $m$ of masses on each loop, at the moment a mass on one loop crosses an intersection with another loop, the corresponding mass on this other loop  crosses the other intersection between the two loops.  Thus no two masses cross an intersection at the same time, and no collisions occur.

Table~\ref{tab:fourier} gives Fourier coefficients found by our action minimization for $m=1,3,5,7$, normalized with $a_1=1$.  As for $m=1$, we checked that the force and acceleration agreed to within $5$ decimal figures at all times.  Movies of these orbits are available at~\cite{web}; a movie of the case $m=1$ also appears at~\cite{ferrarioweb}.  

The special case when $m$ is a multiple of $3$ has some additional symmetries.  By Lagrange symmetry, the set of masses is fixed under translating time by $2\pi/3$.  Therefore, at all times the set of masses is now fixed not just under the four transformations~\eqref{eq:klein}, but also under the cyclic symmetry $x \to y \to z \to x$.  These transformations generate a 12-element symmetry group, consisting of half the orientation-preserving rotations of the cube (formally, a subgroup which is isomorphic to the alternating group $A_4$).  In terms of linear transformations, these are the matrices
\begin{equation} 
\label{eq:a4}
\left\{ 
\mat \pm 1 & & \\ & \pm 1 & \\ & & \pm 1 \rix,  
\mat & \pm 1 & \\ & & \pm 1 \\ \pm 1 & & \rix, 
\mat & & \pm 1 \\ \pm 1 & & \\ & \pm 1 & \rix
\right\} 
\end{equation}
with the restriction that the number of $-1$s in each matrix is even.

This symmetry has amusing consequences.  Recall that the moment of inertia tensor $I$ and the quadrupole moment $Q$ are defined as 
\[ I = \sum_i m_i |\x_i|^2 (\one-\pi_{\x_i}) , \;
Q = \sum_i m_i |\x_i|^2 (3 \pi_{\x_i} - \one) \]
where $\pi_\x$ denotes the projection operator onto the subspace parallel to $\x$, and $\one$ denotes the identity matrix.  Then $I$ and $Q$ must commute with each matrix in the symmetry group~\eqref{eq:a4}.  However, since~\eqref{eq:a4} is an irreducible representation of $A_4$~\cite{FultonH91}, the only such matrices are scalars, i.e., multiples of the identity.  Moreover, since $Q$ has zero trace, it must be identically zero.  Thus, whenever $m$ is a multiple of $3$, $I$ and $Q$ are just as they would be for a sphere, indicating that the mass distribution of these orbits is spherical to second order.  (The moment of inertia is varies with time, however, so these orbits are not counterexamples to Saari's conjecture~\cite{saari} that the only orbits with constant moment of inertia are central configurations.)  

Finally, we note an interesting property of these orbits.  As $m$ grows large, the Fourier coefficients $a_k$, normalized with $a_1=1$, appear to converge (albeit slowly) to asymptotic values for each $k$: Table~\ref{tab:normalized} gives the first 5 nonvanishing Fourier coefficients for $m=61,81,101$, and $121$.  Based on this observation, we conjecture that as $m \to \infty$, these orbits approach a fixed asymptotic shape.  A similar situation seems to hold for the figure-8~\cite{simo1,michael}.

\begin{table}
\[
\begin{array}{|l|l|l|l|l|}
\hline
 k   &  a_k \;(m=1) & a_k \;(m=3) & a_k \;(m=5) & a_k \;(m=7)  \\ \hline 
 1  & 1.00000  & 1.00000 & 1.00000 & 1.00000\\  
 3  & .03282   & -.04629 &-.03991  & -.03335\\
 5  & -.00098  & -.00472 & .00359  &  .00885\\ 
 7  & -.00036  & -.00269 & .00161  &  .00445\\ 
  9 & -.00003  &  .00056 &-.00168  & -.00334\\
 11 &          &  .00010 &-.00043  & -.00039 \\
 13 &          &  .00007 &-.00028  & -.00028\\
 15 &          & -.00002 & .00012  & -.00013 \\ 
 17 &          &         & .00001  & -.00011  \\ 
 19 &          &         & .00001  & -.00008  \\ 
 21 &          &         &         & -.00005  \\ \hline
\end{array} 
\]
\caption{Fourier coefficients $a_k$ for the cubic orbits with $a_1$ normalized to unity.  There are $m$ masses on each loop for a total of $4m$ masses.}
\label{tab:fourier}
\end{table}

\begin{table}
\[
\begin{array}{|l|l|l|l|l|}
\hline
 k   &  a_k \;(m=61) & a_k \;(m=81) & a_k \;(m=101) & a_k \;(m=121) \\ \hline 
1  &1.00000 &1.00000  &1.00000  &1.00000\\
3  &-.02302 &-.02285  &-.02265  &-.02235\\
5  & .01054 & .01025  & .01008  & .00991\\
7  & .00544 & .00522  & .00507  & .00494\\
9  &-.00181 &-.00175  &-.00167  &-.00156\\ \hline
\end{array} 
\]
\caption{Fourier coefficients $a_k$ for the orbits with cubic symmetry
with $a_1$ normalized to one.  Here $n=4m$, with $m$ masses on each loop.}
\label{tab:normalized}
\end{table}

\section{The Criss-cross}

The three-dimensional periodic orbits we discuss above are unstable. In contrast, we  present now a two-dimensional periodic orbit for three masses which we call the ``criss-cross.'' This orbit was  first obtained in 1976 by H\'enon~\cite{henon} as a member of a family of finite angular momentum orbits for equal masses, which are an extension of the zero angular momentum Schubart orbit.  In 1993 it was rediscovered by Moore by searching for orbits with a particular braid type~\cite{moore}, and later shown numerically to exist for different masses by Nauenberg~\cite{michael}. More recently a rigorous proof of existence has been given by Chen~\cite{chen} as part of a family of retrograde orbits.  

H\'enon~\cite{henon} studied the linear stability of the extended Schubart family, and predicted that it contains a finite interval of periodic orbits which are stable in three dimensions.  The extent of this interval, and whether it includes the criss-cross orbit, was not calculated.  However, even if an orbit is stable to first order, it is entirely possible that is not dynamically stable.  Here we provide some numerical evidence that the criss-cross orbit is in fact stable to perturbations, not just in the plane, but in three dimensions.  

The Fourier expansion for the coordinates is given by
\be
x_i(t)=\sum a_{i,k} \cos kt,\; y_i= \sum b_{i,k} \sin kt 
\en
for $i=1,2,3$.
Under the transformation $t \rightarrow t+\pi$, these  coordinates 
change sign, and therefore the sum is restricted to odd values of $k$.
For equal masses, the Fourier coefficients also
satisfy the symmetry condition $a_{2,k}=\pm b_{1,k}$, $b_{2,k}= \pm a_{1,k}$ 
and $b_{3,k}=\pm a_{3,k}$, where the plus sign applies for $k=3,7,11,15\ldots$
and the minus sign applies for $k=1,5,9,13\ldots$  Here masses $1$ and $2$ are 
those orbiting each other in the center, while mass $3$ orbits them both in the retrograde 
direction. 

\begin{figure}
\includegraphics[width=7cm]{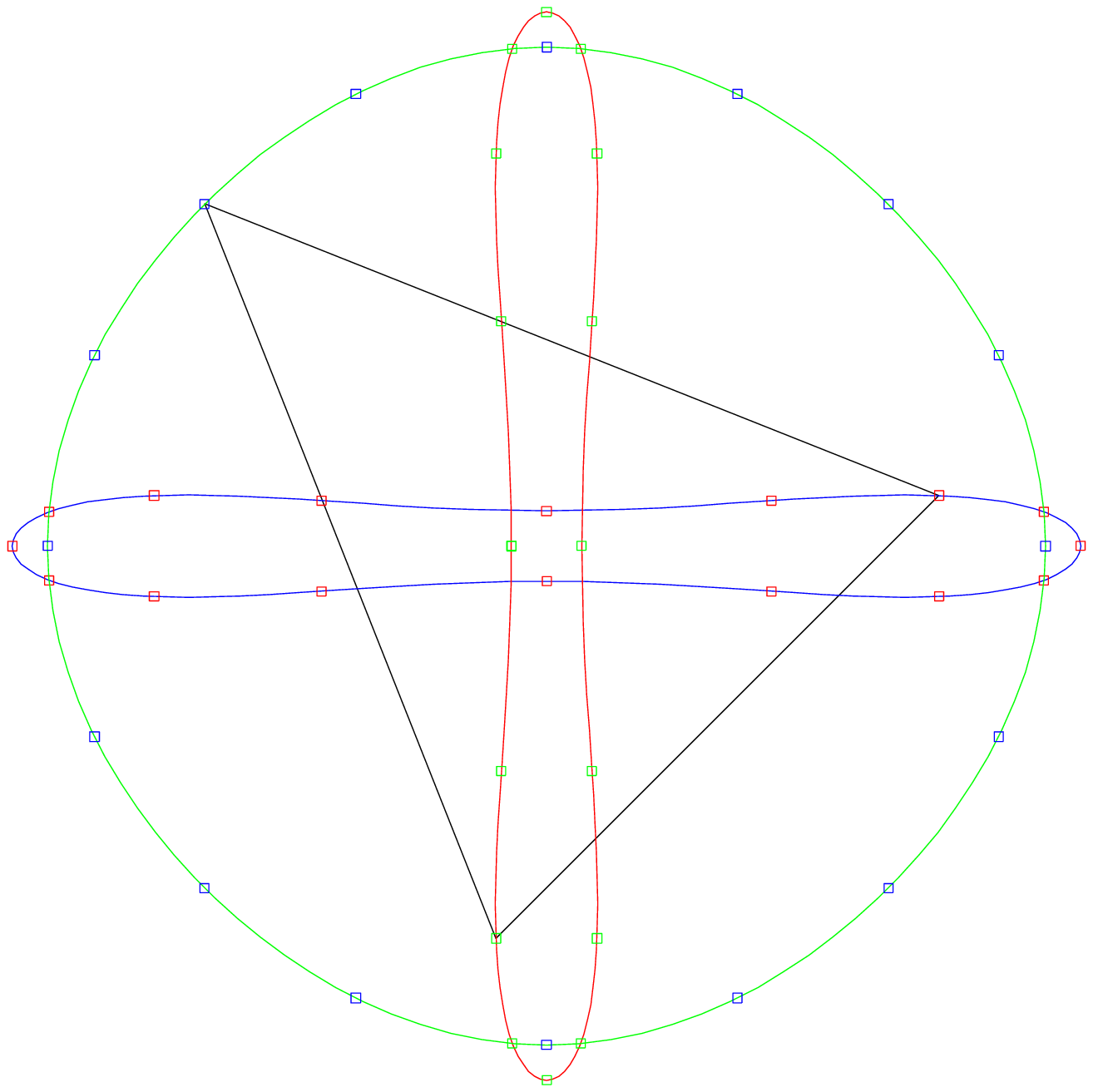}
\includegraphics[width=7cm]{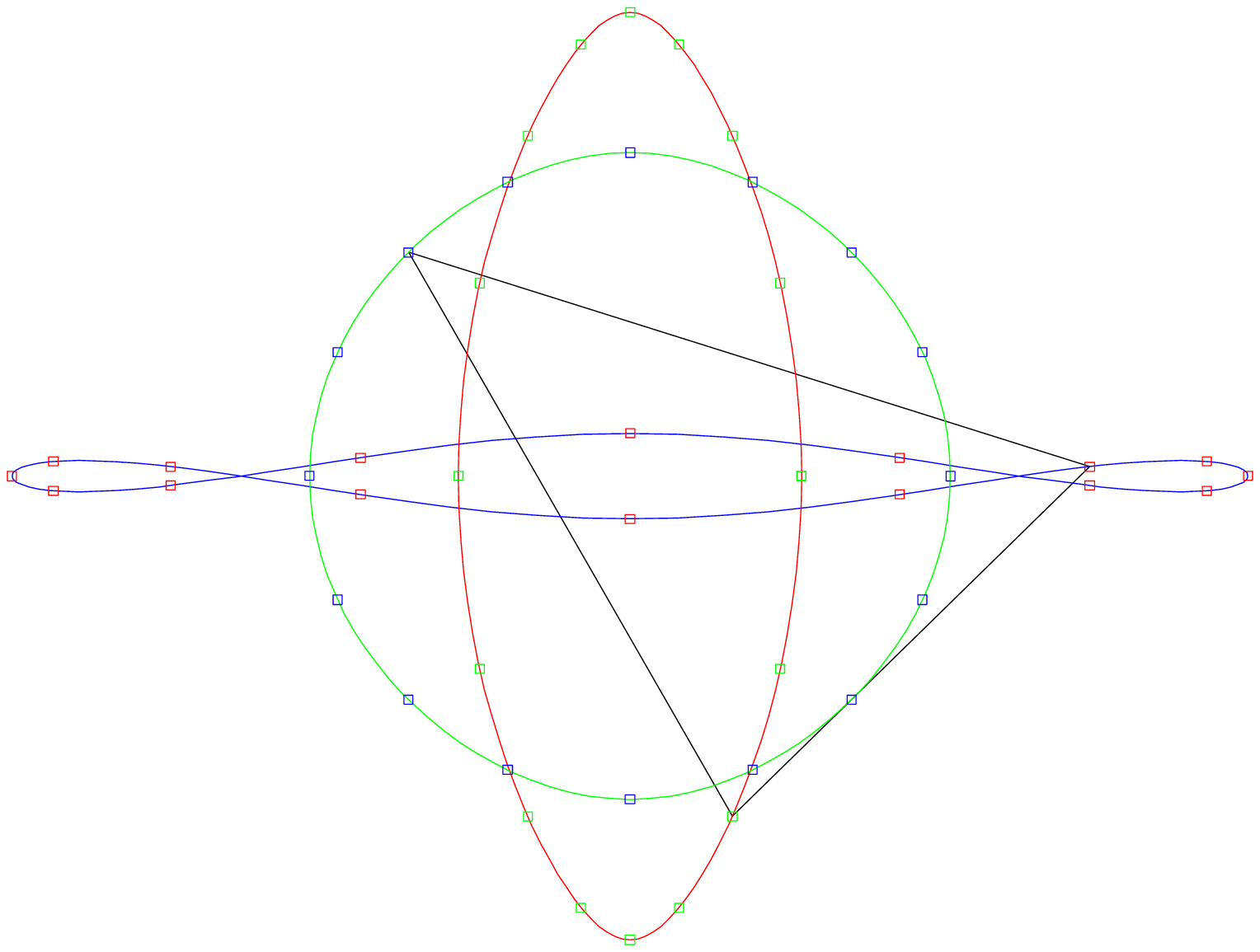}
\caption{
On the left, the criss-cross orbit for three equal masses shown at 16 equal time intervals. Starting with the masses aligned horizontally, at the second time interval they lie at the vertices of an isosceles triangle.  On the right, the corresponding criss-cross orbit with  three masses in the ratio 1:2:3.
}
\label{fig:criss}
\end{figure}

To start the Fourier gradient search of the minimum
of the action integral, we take for the initial Fourier coefficients
$a_{1,1}=1$, $b_{1,1}=0$, and $a_{3,1}=-1$. The resulting
coefficients, which provide Fourier sums that converge to 5 decimal figures,
are shown in Table 3.
In Fig.~\ref{fig:criss} on the left, we show the orbits for three equal masses 
at 16 equal time intervals.  One of the three is nearly circular and rotates in the
opposite direction of the other two orbits.  The initial conditions for the three masses are
\begin{gather*}
 \x_1 =(1.07590,0) , \; \dot\x_1 = (0, 0.19509) \\
 \x_2 =(-0.07095,0) , \; \dot\x_2 = (0,  -1.23187) \\
 \x_3 =(-1.00496,0) , \;  \dot\x_3 = (0, 1.03678)
\end{gather*}
Similar orbits exist for other mass ratios as well.  On the right of Figure~\ref{fig:criss} 
we show the corresponding orbit with three masses in the ratio
1:2:3, showing that similar orbits exists for mass ratios far from $1$.  

\begin{figure}
\includegraphics[width=7cm]{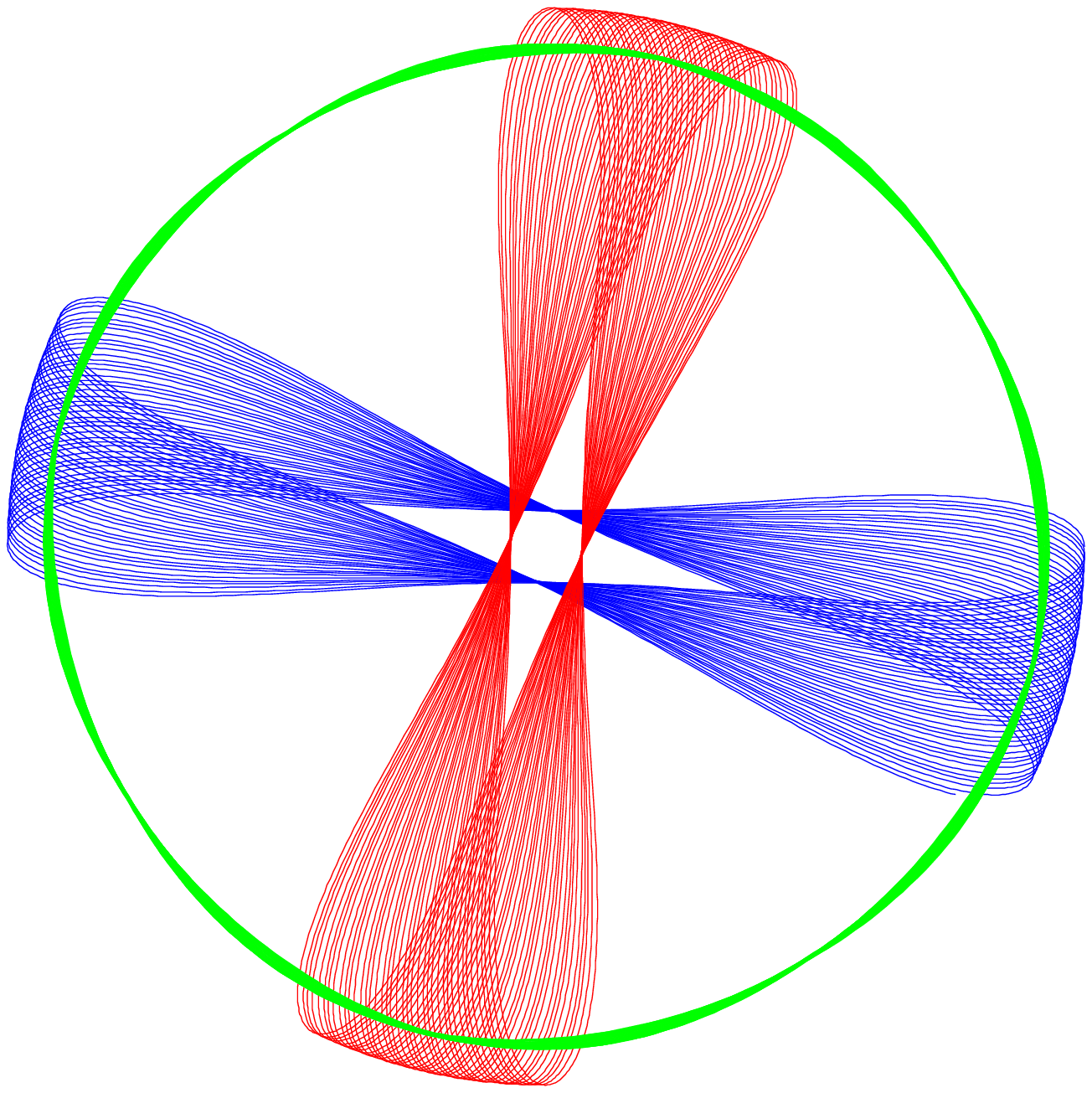}
\includegraphics[width=7cm]{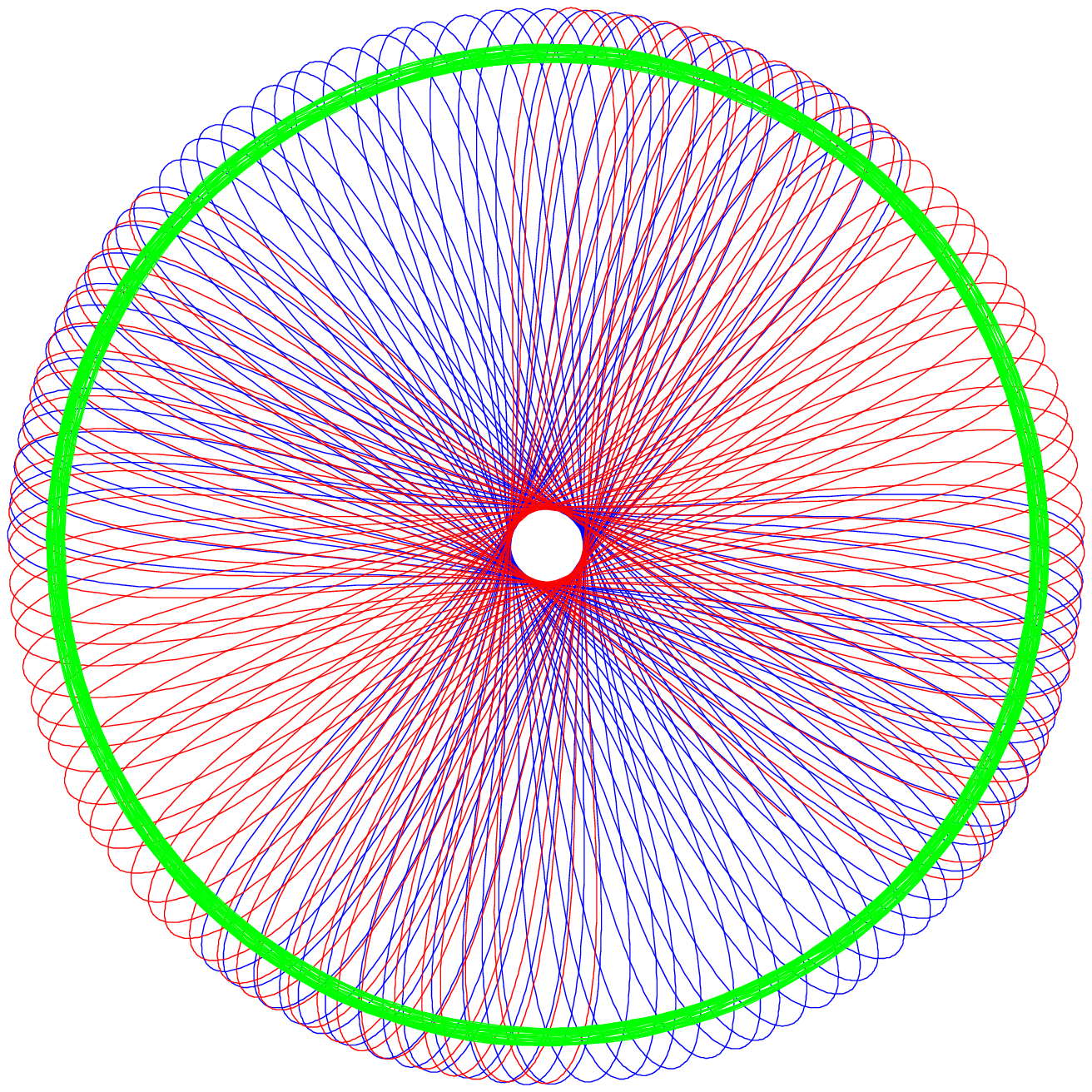}
\caption{On the left, the equal mass criss-cross orbit during 40 periods, for an initial small deviation $\del x=.001$ of the initial position of one of the masses. On the right, the same orbit with a larger perturbation $\del x=.005$}
\label{fig:crissperturbed}
\end{figure}

\begin{figure}
\centerline{\includegraphics[width=8cm]{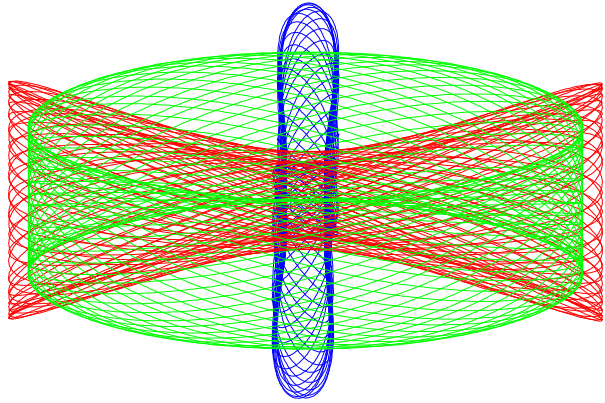}}
\caption{Perturbing the criss-cross along the $z$ direction leads to this beautiful quasiperiodic orbit, indicating that the criss-cross is dynamically stable in three dimensions.  The vertical component of motion has been exaggerated by a factor of 50 for purposes of illustration.}
\label{fig:crissz}
\end{figure}

In Fig.~\ref{fig:crissperturbed} on the left, we show the equal mass criss-cross orbit
during 40 periods for an initial small deviation $\delta x= .001$ of the 
initial position of one of the masses, indicating that the resulting orbit
is quasi-periodic. For longer times the orbit fills out the circle. On the
right side we show a corresponding quasi-periodic orbit which results
for a somewhat larger initial deviation $\delta x= .005 $ of the initial position.  
Similarly, 
a deviation $\delta z=.005$ along the $z$-axis results in the beautiful
quasi-periodic orbit in three dimensions shown edge-wise in Fig.~\ref{fig:crissz}.
To exhibit this behavior we have magnified the axis normal to the plane
of the original periodic orbit by a factor of 50.

These sample results suggest that this orbit is dynamically stable with respect to 
a fairly wide range of perturbations, both in the plane and in three dimensions.  
This stability, and the fact that similar orbits exist for a wide range of mass ratios, suggests that 
triple systems of this type may actually exist in astronomical systems.

\begin{table}
\[
\begin{array}{|l|l|l|l|}
\hline
    k&   a_{1,k} &   b_{1,k}&    a_{3,k} \\ \hline 
    1&   1.09764 &  0.10896 &   -0.98868 \\ 
    3&  -0.02809 &  0.03251 &   -0.00442 \\ 
    5&   0.00724 & -0.00376 &   -0.01100 \\
    7&  -0.00121 &  0.00131 &   -0.00010 \\
    9&   0.00040 & -0.00029 &   -0.00069 \\
   11&  -0.00010 &  0.00010 &   -0.00001 \\
   13&   0.00003 & -0.00003 &   -0.00006 \\
   15&  -0.00001 &  0.00001 &    0.00000 \\ \hline
\end{array}
\]
\caption{Fourier coefficients for the equal mass criss-cross orbit.}
\label{criss-cross}
\end{table}
                                                                                
\end{document}